\documentclass{amsart}
\usepackage{epsfig, amssymb}

\newtheorem{theorem}{Theorem}
\theoremstyle{remark}
\newtheorem*{remark}{Remark}

\newcommand{\R}{{\bf R}} 
\newcommand{\Z}{{\bf Z}} 
\newcommand{\La}{{\mathcal L}}

\title{Lattices in $\R^2$ and finite subsets of a circle}

\author[J.Mostovoy]{Jacob Mostovoy}
\address{Instituto de Matem\'{a}ticas (Unidad Cuernavaca)\\
Universidad Nacional Aut\'{o}noma de M\'{e}xico\\
A.P. 273 Admon. de Correos \# 3\\
C.P. 62251, Cuernavaca, Morelos, MEXICO}
\email{jacob@matcuer.unam.mx}

\date{\today}

\begin{document}

\begin{abstract}
An elementary geometric construction is used to relate the space of 
lattices in $\R^2$ to the space $\exp_3 S^1$ of the subsets of a circle of 
cardinality at most 3. As a consequence we obtain a new proof of known 
statements about the space $\exp_3 S^1$.
\end{abstract}

\maketitle
\subsection*{Spaces of finite subsets}
For $X$  a topological space let $\exp_k X $ be the set of all
finite subsets of $X$ of cardinality at most $k$. There is a map from the 
Cartesian product of $k$ copies of $X$ with itself to $\exp_k X $ 
which sends $(x_1,\ldots, x_k)$ to $\{x_1\}\cup\ldots\cup\{x_k\}$. 
The quotient topology gives $\exp_k X $ the structure of a topological
space. Clearly, $\exp_1 X=X$ for any $X$. 

The simplest non-trivial example is provided by the space 
$\exp_2 S^1$ which is easily seen to be homeomorphic to the M\"{o}bius
band. The space $\exp_3 S^1$ is described by a well-known theorem of R. Bott:
\begin{theorem}{\rm (Bott, \cite{Bott})}\label{thm:Bott}
The space $\exp_3 S^1$ is homeomorphic to a 3-sphere $S^3$.
\end{theorem}
The homeomorphism between $\exp_3 S^1$ and $S^3$ is rather non-obvious. 
This is
illustrated by the following result which is, apparently, due to E. Shchepin.
Consider the canonical embedding $\Delta :S^1\to\exp_3 S^1=S^3$ which sends 
a point $x\in S^1$ to the subset $\{x\}\in\exp_3 S^1$. 
\begin{theorem}{\rm (Shchepin, [unpublished])}\label{thm:Shchepin}
The map $\Delta :S^1\to S^3$ is a trefoil knot.
\end{theorem}

In his proof of Theorem \ref{thm:Bott} Bott used a ``cut-and-paste''
argument. Shchepin's proof of Theorem \ref{thm:Shchepin} was based on a
direct calculation of the fundamental group of 
$S^3\backslash \Delta(S^1)$. The purpose of this note is to show that 
the above theorems can be deduced from well-known facts about lattices
in a plane. 

\subsection*{Spaces of lattices.}

A lattice in $\R^2$ is a subgroup of the vector space $\R^2$ 
generated by two linearly independent vectors.
The set $\La$ of all lattices considered up to multiplication by a non-zero 
real number can be identified with $SL(2,\R)/SL(2,\Z)$ 
and will be considered with this topology.

The space $\La$ can be compactified by adding points which correspond to 
degenerate lattices, that is, subgroups of $\R^2$ generated by one
non-zero vector. Degenerate lattices modulo multiplication 
by a non-zero real number form a circle $S^1$. The topology on the union 
$\widehat{\La}=\La\cup S^1$ can be described as follows. 
Denote by $\widehat{T}$ the space of triangles of perimeter 1 in $\R^2$ 
each of whose angles does not exceed $\pi/2$ and is possibly equal to zero. 
The space $\widehat{T}$ contains degenerate triangles which have one 
side equal to 0 and both angles adjacent to it equal to $\pi/2$. Let 
$T\subset\widehat{T}$ be the subspace of non-degenerate triangles. There is a 
map $p: \widehat{T}\to\widehat{\La}$ which sends a triangle to the lattice 
generated by any two of its sides considered as vectors.
Notice that the restriction of $p$ to $T$ is a continuous map from $T$ onto 
$\La$. We define the topology on $\widehat{\La}$ as the quotient
topology with respect to the map $p$.  

\begin{theorem}\label{thm:lattice} 
The compactification $\widehat{\La}$ of the space of lattices $\La$ is 
homeomorphic to a 3-sphere $S^3$. The inclusion 
$\widehat{\La}\backslash\La\to\widehat{\La}$ is 
a trefoil knot.
\end{theorem} 

The proof of the above theorem (due to D.Quillen) can be found
in \cite{Milnor} on page 84 where a homeomorphism between 
$\La$ and the complement of a trefoil is constructed. This map
can be easily seen to extend to a homeomorphism between $\widehat{\La}$
and $S^3$.

\subsection*{Main theorem} Our main result relates Theorem \ref{thm:lattice} 
to Theorems \ref{thm:Bott} and \ref{thm:Shchepin}.
\begin{theorem}\label{thm:main}
There is a homeomorphism $\Phi : \widehat{\La}\to\exp_3 S^1$
which identifies the circle $\widehat{\La}\backslash\La$ with $\Delta (S^1)$.
\end{theorem} 
\begin{remark}
The circle $S^1$ acts as a subgroup of $PSL(2,\R)$ on  $\La$ (and, in fact, 
on $\widehat{\La}$) rotating the lattices. It will 
be clear from the construction below that $\Phi$ sends the action of $S^1$ on
$\widehat{\La}$ to the action of $SO(2)$ by rotations on $\exp_3 S^1$.
\end{remark}
\begin{proof}
For any non-degenerate lattice $L$ there is a finite number of triangles
$V^{L}_i$ with vertices in $L$ whose sides are generators of $L$ and whose 
angles do not exceed $\pi/2$. There are 12 such triangles for a rectangular 
lattice and 6 for any other lattice, as shown in the figure.
\begin{figure}[ht]
\[\epsffile{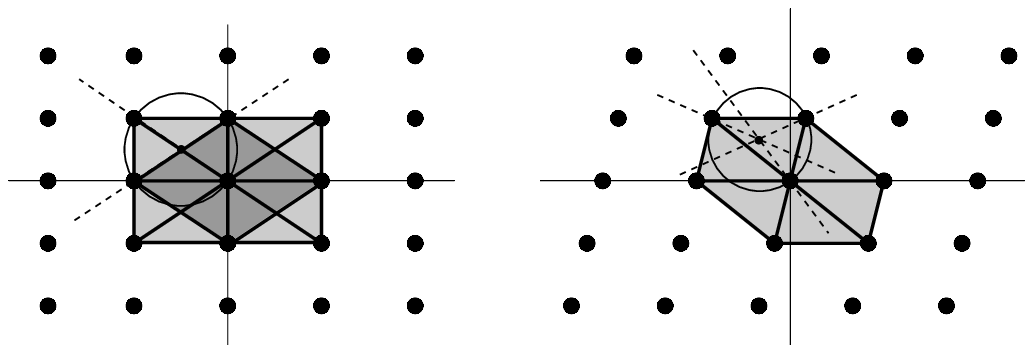}\]
\end{figure}
For each non-degenerate lattice $L$ choose one of the triangles $V^{L}_i$ and
trace 3 lines which connect the vertices of the triangle with the centre 
of its  circumscribed circle. These lines are distinct unless $L$ is 
rectangular in which case two of them coincide. Each line 
determines a point in  ${\bf RP}^1=S^1$ so we obtain a point $\Phi(L)$ in 
$\exp_3 S^1\backslash\Delta (S^1)$ for each $L\in\La$.
It is clear that $\Phi(L)$ depends only on $L$ and not on a particular choice
of a triangle. Now, for a degenerate lattice $L$ define $\Phi(L)$ as the
point of ${\bf RP}^1$ which corresponds to the line 
containing $L$. A straightforward check shows that $\Phi$ is a homeomorphism. 
\end{proof}

\subsection*{Acknowledgments} 
I would like to thank Evgenii Shchepin, Abdelghani
Mouhallil and Alberto Verjovsky for useful discussions.

\end{document}